\documentclass[12pt]{amsart}
\usepackage{amsfonts}
\usepackage{ifthen}
\usepackage{verbatim}
\usepackage{amsthm}
\usepackage{amsmath}
\usepackage{amssymb}
\usepackage{graphicx}
\usepackage{subfigure}
\usepackage{amscd,amssymb,amsthm}
\usepackage{color,xcolor}

\newcounter{minutes}
\divide\time by 60
\newcounter{hours}
\multiply\time by 60 \addtocounter{minutes}{-\time}

\title{Some generalizations of Jensen's inequality}

\author{Slavko Simi\'c}

\address{Mathematical Institute SANU, 11000 Belgrade, Serbia.}\email{ssimic@turing.mi.sanu.ac.rs}


\newtheorem{theorem}[equation]{Theorem}

\newtheorem{lemma}[equation]{Lemma}

\newtheorem{corollary}[equation]{Corollary}
\newtheorem{remark}[equation]{Remark}

\newcommand{\beq}{\begin{equation}}
\newcommand{\eeq}{\end{equation}}

\numberwithin{equation}{section}

\begin{document}

\def\thefootnote{}
\footnotetext{ \texttt{\tiny File:~\jobname .tex,
          printed: \number\year-\number\month-\number\day,
          \thehours.\ifnum\theminutes<10{0}\fi\theminutes}
} \makeatletter\def\thefootnote{\@arabic\c@footnote}\makeatother

\begin{abstract}
In this article we give some improvements and generalizations of
the famous Jensen's and Jensen-Mercer inequalities for twice
differentiable functions, where convexity property of the target
function is not assumed in advance. They represents a refinement
of these inequalities in the case of convex/concave functions with
numerous applications in Theory of Means and Probability and
Statistics.
\end{abstract}

\maketitle
\section{Introduction}

\bigskip

Recall that the Jensen functional $J_n(\bold p, \bold x; f)$ is
defined on an interval $ I\subseteq \mathbb R$ by
$$
J_n(\bold p,\bold x;f):=\sum_1^n p_i f(x_i)-f(\sum_1^n p_i x_i),
$$
where $f: I\to \mathbb R$, \ $\bold x=(x_1, x_2,\cdots, x_n)\in
I^n$ and $\bold p=\{p_i\}_1^n $ is a positive weight sequence.

\bigskip

Another well known assertion is the famous

\bigskip

{\bf  Jensen's inequality} (\cite{hlp}) \ {\it If $f$ is twice
continuously differentiable function and $f''\ge 0$ on an interval
$I$, then $f$ is convex on $I$ and the inequality
$$
0\le J_n(\bold p,\bold x;f)
$$
holds for each $\bold x:=(x_1, ..., x_n)\in I^n$ and any positive
weight sequence $\bold p:=\{p_i\}_1^n$ with $\sum_1^n p_i=1$.

If $f''\le 0$ on $I$, then $f$ is a concave function on $I$ and}
$$
J_n(\bold p,\bold x;f)\le 0.
$$

\bigskip

Its counterpart is given by the following

\bigskip

{\bf Jensen-Mercer inequality} (\cite{m}) \ {\it Let $\phi:
[a,b]\subseteq \mathbb R \to \mathbb R$ be a convex function and
$x_i\in [a,b], i=1,2,...,n$. Then}

\begin{equation}\label{eq1}
\phi(a+b-\sum_1^n p_ix_i)\le \phi(a)+\phi(b)-\sum_1^n
p_i\phi(x_i).
\end{equation}

\bigskip

Our first task in this paper is to find some global upper bounds
for these inequalities. We prove the following

\bigskip

{\it Let $f$ be a convex function on an interval $I$ and $x_i\in
[a,b]\subset I$. Then}

$$
0\le \sum_1^n p_i f(x_i)-f(\sum_1^n p_i x_i)\le
f(a)+f(b)-2f(\frac{a+b}{2});
$$

\bigskip

$$
0\le f(a)+f(b)-\sum_1^n p_if(x_i)-f(a+b-\sum_1^n p_ix_i)\le
2(f(a)+f(b)-2f(\frac{a+b}{2})).
$$

\bigskip

Those bounds can be improved by the {\it characteristic} number
$c(f)$ of the convex function $f$ (\cite {ss}), to the next

$$
0\le \sum_1^n p_i f(x_i)-f(\sum_1^n p_i x_i)\le
c(f)[f(a)+f(b)-2f(\frac{a+b}{2})];
$$

\bigskip

$$
0\le f(a)+f(b)-\sum_1^n p_if(x_i)-f(a+b-\sum_1^n p_ix_i)\le
(1+c(f))[f(a)+f(b)-2f(\frac{a+b}{2})].
$$

\bigskip

As an example, we shall calculate characteristic number for the
power function:

$$
c(x^s)=\begin{cases} 1,& s < 0; \\
                    (1-s)s^{s/(1-s)}/(2^{1-s}-1), & 0<s<1; \\
                    (s-1)s^{s/(1-s)}/(1-2^{1-s}), & s>1. \\
                    \end{cases}
$$

\bigskip

Our second main task is to investigate the possibility of a form
of Jensen's and Jensen-Mercer inequalities for functions which are
not necessarily convex/concave on $I$.

The sole condition will be that the second derivative of the
target function exists locally i.e., on a closed interval
$E:=[a,b]\subset I$. Since it is continuous on a closed interval,
there exist numbers $m_f(E)=m(a,b;f):=\min_{t\in E}f''(t)$ and
$M_f(E)=M(a,b;f):=\max_{t\in E}f''(t)$. Those numbers will play an
important role in the sequel.

\bigskip

For instance, let $f\in C^{(2)}(E)$ and $x_i\in E, \ i=1,2,...,n$.
Then

$$
\frac{1}{2}m_f(E) J_n(\bold p,\bold x;x^2)\le J_n(\bold p,\bold
x;f)\le \frac{1}{2}M_f(E) J_n(\bold p,\bold x;x^2).
$$

\bigskip

Note that this inequality represents an improvement of Jensen's
inequality for convex functions since in this case we have $0\le
m_f(E)<M_f(E)$.

\bigskip

\section{Results and Proofs}

\bigskip

We firstly determine some global upper bounds for Jensen's and
Jensen-Mercer inequalities.

\bigskip

\begin{theorem}\label{thm1} \ Let $f$ be a convex function on $I$
and $\bold x\in [a,b]\subset I$. Then

\begin{equation}\label{eq2}
0\le J_n(\bold p,\bold x;f)=\sum_1^n p_i f(x_i)-f(\sum_1^n p_i
x_i)\le f(a)+f(b)-2f(\frac{a+b}{2});
\end{equation}

\bigskip

\begin{equation}\label{eq3}
0\le f(a)+f(b)-\sum_1^n p_if(x_i)-f(a+b-\sum_1^n p_ix_i)\le
2[f(a)+f(b)-2f(\frac{a+b}{2})],
\end{equation}

independently of $\bold p$.

\end{theorem}

\bigskip

\begin{proof} \ We obtain a simple proof of (\ref{eq2}) directly
from Jensen-Mercer inequality and

\begin{lemma}\label{l1} \cite{s} \ Let $h$ be a convex function on $E=[a,b]$ and, for some $x,y\in
E, \ x+y=a+b$, then

$$
2h(\frac{a+b}{2})\le h(x)+h(y)\le h(a)+h(b).
$$
\end{lemma}

\bigskip

Namely, writing this inequality in the form

$$
\sum_1^n p_i f(x_i)-f(\sum_1^n p_i x_i)\le f(a)+f(b)-(f(\sum_1^n
p_ix_i)+f(a+b-\sum_1^n p_ix_i)),
$$

the proof follows by Lemma \ref{l1}.

\bigskip

For the proof of the assertion (\ref{eq3}), note that if $x_i\in
[a,b]$ then also $y_i:=a+b-x_i\in [a,b]$. Hence, by (\ref{eq2})
and Lemma \ref{l1}, we get

$$
f(a)+f(b)-2f(\frac{a+b}{2})\ge \sum_1^n p_i f(y_i)-f(\sum_1^n p_i
y_i)
$$

$$
=\sum_1^n p_i f(a+b-x_i)-f(\sum_1^n p_i(a+b- x_i))\ge \sum_1^n p_i
[2f(\frac{a+b}{2})-f(x_i)]-f(\sum_1^n (a+b-p_i x_i))
$$

$$
=f(a)+f(b)-\sum_1^n p_if(x_i)-f(a+b-\sum_1^n
p_ix_i)-[f(a)+f(b)-2f(\frac{a+b}{2})],
$$

and the proof is done.

\end{proof}

\bigskip

Those bounds can be improved by the {\it characteristic} number
$c(f)$ of the convex function $f$ (\cite {ss}), to the following

\begin{equation}\label{eq4}
0\le \sum_1^n p_i f(x_i)-f(\sum_1^n p_i x_i)\le
c(f)[f(a)+f(b)-2f(\frac{a+b}{2})];
\end{equation}

\bigskip

\begin{equation}\label{eq5}
0\le f(a)+f(b)-\sum_1^n p_if(x_i)-f(a+b-\sum_1^n p_ix_i)\le
(1+c(f))[f(a)+f(b)-2f(\frac{a+b}{2})],
\end{equation}

where the characteristic number $c(f), c(f)\in [1/2,1]$ is defined
by

$$
c(f):=\sup_{p,q;a,b}\frac{pf(a)+qf(b)-f(pa+qb)}{f(a)+f(b)-2f(\frac{a+b}{2})}.
$$

\bigskip

To find $c(f)$ for a concrete convex function $f$ is not an easy
task. Nevertheless, by direct calculation we obtain

$$
c(x^2)=\sup_{p,q} 2pq=1/2.
$$

We shall determine now the value of this constant for some classes
of functions.

\bigskip

For this cause, recall the definitions of slowly varying and
rapidly varying functions (cf. [BGT]).

\bigskip

 {\bf Definition}\ Let the function $f$ be defined on $I:=[a,+\infty)$.

 \bigskip

It is said that $f$ is slowly varying if
$\lim_{x\to\infty}\frac{f(tx)}{f(x)}=1$ for any $t>0$.

\bigskip

If $\lim_{x\to\infty}\frac{f(tx)}{f(x)}=\infty$ for any $t>1$,
then $f$ is a rapidly varying function.

\bigskip

\begin{theorem} \ Let $f(a+x):=g_a(x)$ be a slowly or rapidly
varying function. Then $c(f)=1$.

\begin{proof} \ Denote

$$
H:=\frac{pf(a)+qf(b)-f(pa+qb)}{f(a)+f(b)-2f(\frac{a+b}{2})}=\frac{pf(a)+qg_a(x)-g_a(qx)}{f(a)+g_a(x)-2g_a(\frac{x}{2})},
$$

with $x=b-a$.

\bigskip

Since $f$ is a convex function, so is $g_a(x)$.
 Hence $\lim_{x\to\infty}g_a(x)$ can be $0,c$ or $\pm\infty$.

 \bigskip

 In the first two cases we obtain at once that $\lim_{x\to\infty}H=p$. Since
 $g_a(x)$ is also slowly varying, in the third case we get

 $$
\lim_{x\to\infty}H=\frac{pf(a)/g_a(x)+q-g_a(qx)/g_a(x)}{f(a)/g_a(x)+1-2g_a(\frac{x}{2})/g_a(x)}=\frac{q-1}{-1}=p.
 $$

 \bigskip

 As concerns the class of rapidly varying functions, note that $\lim_{x\to\infty}\frac{f(tx)}{f(x)}=0$
 for $0<t<1$, which can be easily proven by the change of variable
 $tx\to x, 1/t\to t$.

 Therefore, in this case we have

 $$
\lim_{x\to\infty}H=\frac{pf(a)/g_a(x)+q-g_a(qx)/g_a(x)}{f(a)/g_a(x)+1-2g_a(\frac{x}{2})/g_a(x)}=q.
 $$

 Since $p$ and $q$ are arbitrary weights, we conclude that $c(f)=1$
 in both cases. For instance,

 $$
c(-\log x)=c(e^{-x})=c(e^x)=c(x^x)=1.
 $$
\end{proof}

\end{theorem}

\bigskip

Our next contribution is an evaluation of the characteristic
number for the power function.

\bigskip

\begin{theorem}\label{thm7} \ We have

$$
c(x^s)=\begin{cases} 1,& s < 0; \\
                    (1-s)s^{s/(1-s)}/(2^{1-s}-1), & 0<s<1; \\
                    (s-1)s^{s/(1-s)}/(1-2^{1-s}), & s>1. \\
                    \end{cases}
$$

\end{theorem}

\begin{proof} \ Main tool for the proof of this and similar theorems will be
the following useful assertion.

\bigskip

\begin{lemma} \label{l2}{\rm \cite[Theorem 1.25]{avv}}.
For $-\infty<a<b<\infty$, let $f,\,g: [a,b]\rightarrow \mathbb{R}$
be continuous on $[a,b]$, and be differentiable on $(a,b)$, and
let $g'(x)\neq 0$ on $(a,b)$. If $f'(x)/g'(x)$ is
increasing(deceasing) on $(a,b)$, then so are
\begin{eqnarray*}
\frac{f(x)-f(a)}{g(x)-g(a)}\,\,\,\,\,\,\,and\,\,\,\,\,\,\,\,\frac{f(b)-f(x)}{g(b)-g(x)}.
\end{eqnarray*}
If $f'(x)/g'(x)$ is strictly monotone, then the monotonicity in
the conclusion is also strict.
\end{lemma}

\bigskip

Let $a,p,q\in\mathbb R^+, \ p+q=1,\ p\neq q; \ x\in (a, +\infty) $
and $s\in (0,1)\cup(1,2)\cup(2,+\infty)$.

\bigskip

Denote $f_1(x)=(q+pa/x)^{s-1}; \ g_1(x)=((1+a/x)/2)^{s-1}$.

Since

$$
\frac{f'_1(x)}{g'_1(x)}=2p\frac{(q+pa/x)^{s-2}}{((1+a/x)/2)^{s-2}}=2p\Bigl(\frac{pa+qx}{(a+x)/2}\Bigr)^{s-2},
$$

by Lemma (\ref{l2}) we conclude that the expression

$$
\frac{f_1(x)-f_1(a)}{g_1(x)-g_1(a)}=\frac{(q+pa/x)^{s-1}-1}{((1+a/x)/2)^{s-1}-1}=\frac{x^{s-1}-(pa+qx)^{s-1}}{x^{s-1}-((a+x)/2)^{s-1}},
$$

\bigskip

is monotone increasing for $q>p, \ s\in (2,+\infty)$ or $p>q, \
s\in (0,1)\cup (1,2)$ and monotone decreasing otherwise.

\bigskip

Denote now $f_2(x)=qx^s-(pa+qx)^s; \ g_2(x)=x^s-2((a+x)/2)^s$.

Since

$$
\frac{f'_2(x)}{g'_2(x)}=q\frac{x^{s-1}-(pa+qx)^{s-1}}{x^{s-1}-((a+x)/2)^{s-1}},
$$

we conclude the same for

$$
\frac{f_2(x)-f_2(a)}{g_2(x)-g_2(a)}=\frac{pa^s+qx^s-(pa+qx)^s}{a^s+x^s-2((a+x)/2)^s}:=H(x).
$$

Hence, the maximum of $H(x)$ is attained at the endpoints of
$(a,+\infty)$.

\bigskip

We have

$$
\lim_{x\to a}H(x)=2pq; \ \lim_{x\to
+\infty}H(x)=\frac{q-q^s}{1-2^{1-s}}.
$$

\bigskip

Because $\max_q(2pq)=1/2$ is the least possible value of $c(f)$,
we see that

$$
c(x^s)=\max_q(q-q^s)/(1-2^{1-s}),
$$

and the proof follows.

\bigskip

For $x\in (0,b)$, putting

$$
f_1(x)=(p+qb/x)^{s-1}, \ g_1(x)=((1+b/x)/2)^{s-1};
$$
$$
f_2(x)=px^s-(px+qb)^s; \ g_2(x)=x^s-2((b+x)/2)^s,
$$

and repeating the above procedure, we obtain the same result.

\bigskip

If $s<0$, we have $\lim_{x\to\infty}x^s=0$. Hence $c(x^s)=1$
according to the previous theorem.

\end{proof}

\bigskip

\begin{remark} \ The described method can be applied for
evaluation of the characteristic number of other convex functions.

For example, it can be proved that $c(x\log x)=(e\log 2)^{-1}$.

\end{remark}

\bigskip

Our next achievement is the form of Jensen's and Jensen-Mercer
inequalities for non-convex functions.

\bigskip

\begin{theorem}\label{thm3} \ Let $g\in C^{(2)}(E)$ and $\bold x\in E:=[a,b]\subset\mathbb R$.

\bigskip

Then

$$
\frac{1}{2}m_f(E) J_n(\bold p,\bold x;x^2)\le J_n(\bold p,\bold
x;g)\le \frac{1}{2}M_f(E) J_n(\bold p,\bold x;x^2).
$$

\bigskip

where $m_f(E):=\min_{t\in E}g''(t)$ and $M_f(E):=\max_{t\in
E}g''(t)$.
\end{theorem}

\bigskip

\begin{proof} \ For a given $g\in C^{(2)}(E)$, define an auxiliary
function $f$ by $f(x):= g(x)-m_g(E) x^2/2$. Since
$f''(x)=g''(x)-m_g(E)\ge 0$, we see that $f$ is a convex function
on $E$. Therefore, applying Jensen's inequality, we obtain

$$
0\le J_n(\bold p,\bold x;f)=J_n(\bold p,\bold
x;g)-\frac{1}{2}m_g(E)J_n(\bold p,\bold x;x^2).
$$

\bigskip

On the other hand, taking the auxiliary function $f$ as
$f(x)=M_g(E)x^2/2-g(x)$, we see that it is also convex on $E$.

Applying Jensen's inequality again, we get

$$
0\le J_n(\bold p,\bold x;f)=\frac{1}{2}M_g(E)J_n(\bold p,\bold
x;x^2)-J_n(\bold p,\bold x;g),
$$

and the proof is done.

\end{proof}

\bigskip

Another form is possible.

\bigskip

\begin{theorem}\label{thm4} \ Let $g\in C^{(2)}(E)$ and $\bold x\in E:=[a,b]\subset\mathbb R$.

\bigskip

Then

$$
g(a)+g(b)-2g(\frac{a+b}{2})+\frac{1}{4}M_g(E) [2J_n(\bold p,\bold
x;x^2)-(b-a)^2]
$$
$$
\le J_n(\bold p,\bold x;g)\le
$$
$$
g(a)+g(b)-2g(\frac{a+b}{2})+\frac{1}{4}m_g(E) [2J_n(\bold p,\bold
x;x^2)-(b-a)^2].
$$
\end{theorem}

\bigskip

\begin{proof} \ Applying the same auxiliary functions to the
converse of Jensen's inequality (\ref{eq2}), we obtain the desired
result.

\end{proof}

\bigskip

Two-sided improvement of Jensen's inequality is given by the next

\begin{theorem} \ Let $f\in C^{(2)}(E)$ be a convex function and $\bold x\in E:=[a,b]\subset\mathbb R$.

\bigskip

Then

$$
\frac{m_f(E)}{m_f(E)+M_f(E)}[f(a)+f(b)-2f(\frac{a+b}{2})]+\frac{m_f(E)M_f(E)}{m_f(E)+M_f(E)}(J_n(\bold
p,\bold x;x^2)-\frac{1}{4}(b-a)^2)
$$
$$
\le J_n(\bold p,\bold x;f)\le
$$
$$
\frac{M_f(E)}{m_f(E)+M_f(E)}[f(a)+f(b)-2f(\frac{a+b}{2})]+\frac{m_f(E)M_f(E)}{m_f(E)+M_f(E)}(J_n(\bold
p,\bold x;x^2)-\frac{1}{4}(b-a)^2).
$$
\end{theorem}

\bigskip

\begin{proof} \ Adjusting the right-hand parts of Theorem
\ref{thm3} and Theorem \ref{thm4}, we obtain

$$
 J_n(\bold p,\bold x;f)\le
\frac{M_f(E)}{m_f(E)+M_f(E)}[f(a)+f(b)-2f(\frac{a+b}{2})+\frac{1}{4}m_f(E)
[2J_n(\bold p,\bold x;x^2)-(b-a)^2]]
$$
$$
+\frac{m_f(E)}{m_f(E)+M_f(E)}[\frac{1}{2}M_f(E)J_n(\bold p,\bold
x;x^2)]
$$
$$
=\frac{M_f(E)}{m_f(E)+M_f(E)}[f(a)+f(b)-2f(\frac{a+b}{2})]+\frac{m_f(E)M_f(E)}{m_f(E)+M_f(E)}(J_n(\bold
p,\bold x;x^2)-\frac{1}{4}(b-a)^2).
$$

\bigskip

Similarly, adjusting left-hand sides we get

$$
 J_n(\bold p,\bold x;f)\ge
\frac{m_f(E)}{m_f(E)+M_f(E)}[f(a)+f(b)-2f(\frac{a+b}{2})+\frac{1}{4}M_f(E)
[2J_n(\bold p,\bold x;x^2)-(b-a)^2]]
$$
$$
+\frac{M_f(E)}{m_f(E)+M_f(E)}[\frac{1}{2}m_f(E)J_n(\bold p,\bold
x;x^2)]
$$
$$
=\frac{m_f(E)}{m_f(E)+M_f(E)}[f(a)+f(b)-2f(\frac{a+b}{2})]+\frac{m_f(E)M_f(E)}{m_f(E)+M_f(E)}(J_n(\bold
p,\bold x;x^2)-\frac{1}{4}(b-a)^2),
$$

and the proof follows.

\bigskip
\end{proof}

\bigskip

A simple consequence of the previous theorem is another converse
of Jensen's inequality.

\bigskip

\begin{corollary} \ Because $J_n(\bold p,\bold x;x^2)\le
\frac{1}{4}(b-a)^2$, we obtain

\bigskip

\begin{equation}
 J_n(\bold p,\bold x;f)\le \frac{M_f(E)}{m_f(E)+M_f(E)}[f(a)+f(b)-2f(\frac{a+b}{2})],
\end{equation}

\end{corollary}

\bigskip

\begin{remark} \ Since $\frac{M_f(E)}{m_f(E)+M_f(E)}\in
[\frac{1}{2}, 1]$, it is interesting to compare this result with
(\ref{eq4}).
\end{remark}

\bigskip

A non-convex variant of the Jensen-Mercer inequality follows.

\bigskip

\begin{theorem}\label{thm6} \ Let $g\in C^{(2)}(E)$ and $\bold x\in E:=[a,b]\subset\mathbb R$.

\bigskip

Then

$$
\frac{1}{2}m_g(E)[2(\sum_1^n p_ix_i-a)(b-\sum_1^n
p_ix_i)-J_n(\bold p,\bold x;x^2)]
$$
$$
\le g(a)+g(b)-\sum_1^n p_ig(x_i)-g(a+b-\sum_1^n p_ix_i)\le
$$
$$
\frac{1}{2}M_g(E)[2(\sum_1^n p_ix_i-a)(b-\sum_1^n
p_ix_i)-J_n(\bold p,\bold x;x^2)].
$$

\end{theorem}

\begin{proof}\ Applying Jensen-Mercer inequality

$$
0\le f(a)+f(b)-\sum_1^n p_if(x_i)-f(a+b-\sum_1^n
p_ix_i):=K_n(\bold p,\bold x;f)
$$

to the convex function $f(x)=g(x)-\frac{1}{2}m_g(E)x^2$, we get

$$
0\le K_n(\bold p,\bold x;g)-\frac{1}{2}m_g(E)K_n(\bold p,\bold
x;x^2)
$$
$$
=K_n(\bold p,\bold x;g)-\frac{1}{2}m_g(E)[a^2+b^2-(a+b-\sum_1^n
p_ix_i)^2-\sum_1^n p_ix_i^2]
$$
$$
=K_n(\bold p,\bold x;g)-\frac{1}{2}m_g(E)[-2ab+2(a+b)\sum_1^n
p_ix_i)-2(\sum_1^n p_ix_i)^2-(\sum_1^n p_ix_i^2-(\sum_1^n
p_ix_i)^2)]
$$
$$
=K_n(\bold p,\bold x;g)-\frac{1}{2}m_g(E)[ 2(\sum_1^n
p_ix_i-a)(b-\sum_1^n p_ix_i)-J_n(\bold p,\bold x;x^2)].
$$

\bigskip

Consequently, for the function $f(x)=\frac{1}{2}M_g(E)x^2-g(x)$ we
obtain

$$
0\le \frac{1}{2}M_g(E)[ 2(\sum_1^n p_ix_i-a)(b-\sum_1^n
p_ix_i)-J_n(\bold p,\bold x;x^2)]-K_n(\bold p,\bold x;g),
$$

and the proof is done.
\end{proof}

\bigskip

 \section{Applications}

\bigskip

{\bf General means} \ Most known general means are

$$
\mathcal{A}(\bold w,\bold x):=\sum w_i x_i;
$$
$$
\mathcal{G}(\bold w,\bold x):=\prod x_i^{w_i};
$$
$$
\mathcal{H}(\bold w,\bold x):=(\sum w_i/x_i)^{-1},
$$

i.e., arithmetic, geometric and harmonic mean, respectively.

\bigskip

Here $\bold x=\{x_i\}_1^n$ denotes an arbitrary sequence of
positive numbers and $\bold w=\{w_i\}_1^n$ is a corresponding
weight sequence.

\bigskip

The famous $\mathcal{A}-\mathcal{G}-\mathcal{H}$ inequality says
that

$$
0\le \mathcal{H}(\bold w,\bold x)\le \mathcal{G}(\bold w,\bold
x)\le \mathcal{A}(\bold w,\bold x).
$$

\bigskip

It is proved in \cite{sim} that $1\le\mathcal{A}/\mathcal{H}\le
(a+b)^2/4ab$, whenever $\bold x\in [a,b]\subset \mathbb R^+$.

\bigskip

The same bounds hold for other
$\mathcal{A}-\mathcal{G}-\mathcal{H}$ quotients.

\bigskip

\begin{theorem} \ Let $\bold x\in [a,b]\subset \mathbb R^+$. Then

$$
1\le\frac{\mathcal{A}(\bold w,\bold x)}{\mathcal{H}(\bold w,\bold
x)}\le \frac{(a+b)^2}{4ab};
$$
$$
1\le\frac{\mathcal{A}(\bold w,\bold x)}{\mathcal{G}(\bold w,\bold
x)}\le \frac{(a+b)^2}{4ab};
$$
$$
1\le\frac{\mathcal{G}(\bold w,\bold x)}{\mathcal{H}(\bold w,\bold
x)}\le \frac{(a+b)^2}{4ab};
$$
\end{theorem}

\begin{proof} \ Since $f(x)=-\log x$ is a convex function on $\mathbb
R^+$, using Theorem \ref{thm1} we get

$$
\log(\sum w_i x_i)-\sum w_i\log x_i\le 2\log \frac{a+b}{2}-\log
a-\log b,
$$

that is,

$$
\log [\frac{\mathcal{A}(\bold w,\bold x)}{\mathcal{G}(\bold
w,\bold x)}]\le \log [\frac{(a+b)^2}{4ab}],
$$

and the proof follows.

\vskip 1cm

Finally,

$$
1\le \frac{\mathcal{G}(\bold w,\bold x)}{\mathcal{H}(\bold w,\bold
x)}=\frac{\mathcal{A}(\bold w,\bold x)}{\mathcal{H}(\bold w,\bold
x)}/\frac{\mathcal{A}(\bold w,\bold x)}{\mathcal{G}(\bold w,\bold
x)}\le \frac{(a+b)^2}{4ab}.
$$

\end{proof}

\bigskip

Similar converses are valid for the
$\mathcal{A}-\mathcal{G}-\mathcal{H}$ differences.

\bigskip

\begin{theorem}  \ Let $\bold x\in [a,b]\subset \mathbb R^+$. Then

$$
0\le {\mathcal{A}(\bold w,\bold x)}-{\mathcal{G}(\bold w,\bold
x)}\le (\sqrt b-\sqrt a)^2;
$$
$$
0\le {\mathcal{A}(\bold w,\bold x)}-{\mathcal{H}(\bold w,\bold
x)}\le (\sqrt b-\sqrt a)^2;
$$
$$
0\le {\mathcal{G}(\bold w,\bold x)}-{\mathcal{H}(\bold w,\bold
x)}\le (\sqrt b-\sqrt a)^2.
$$
\end{theorem}

\bigskip

For example, taking $f(x)=e^x$ and applying Theorem \ref{thm1}, we
obtain

$$
\sum w_i e^{x_i}-e^{\sum w_i x_i}\le e^a +
e^b-2e^{\frac{a+b}{2}}=(e^{b/2}-e^{a/2})^2.
$$

\bigskip

Now, change of variable $\bold x\to\log \bold x; a\to\log a,
b\to\log b$ gives the result.

Rest of the proof is left to the reader.

\bigskip

Notion of $\mathcal{A}-\mathcal{G}-\mathcal{H}$ means is
generalized by the power mean $\mathcal{P}_{\alpha}$ of order
$\alpha\in\mathbb R$, defined as

$$
\mathcal{P}_{\alpha}(\bold x,\bold w):=(\sum w_i
x_i^{\alpha})^{1/{\alpha}}.
$$

\bigskip

Hence,

$$
\mathcal{P}_{-1}(\bold x,\bold w)=\mathcal{H}(\bold x,\bold w), \
\mathcal{P}_1(\bold x,\bold w)=\mathcal{A}(\bold x,\bold w),
$$

and

$$
\mathcal{P}_0(\bold x,\bold w)=\lim_{{\alpha}\to
0}\mathcal{P}_{\alpha}(\bold x,\bold w)=\mathcal{G}(\bold x,\bold
w).
$$

\bigskip

It is well known (\cite{hlp}) that power means are monotone
increasing in $\alpha$.

\bigskip

We give now an estimation of a difference of power means.

\bigskip

\begin{theorem} \ For $0<\alpha<1$ and $\bold x\in [a,b]$, we have

\bigskip

\begin{equation}\label{eq6}
 0\le \mathcal{A}(\bold x,\bold w)-\mathcal{P}_{\alpha}(\bold
x,\bold w)\le
2(1-\alpha)\alpha^{\frac{\alpha}{1-\alpha}}/(1-2^{\frac{\alpha
-1}{\alpha}})[\frac{a+b}{2}-\Bigl(\frac{a^{\alpha}+b^{\alpha}}{2}\Bigr)^{1/\alpha}].
\end{equation}

\bigskip

For $\alpha>1$, we have

\begin{equation}\label{eq7}
 0\le \mathcal{P}_{\alpha}(\bold
x,\bold w)-\mathcal{A}(\bold x,\bold w)\le
2(\alpha-1)\alpha^{\frac{\alpha}{1-\alpha}}/(2^{\frac{\alpha
-1}{\alpha}}-1)[\Bigl(\frac{a^{\alpha}+b^{\alpha}}{2}\Bigr)^{1/\alpha}-\frac{a+b}{2}].
\end{equation}

\end{theorem}

\begin{proof} \ By Theorem \ref{thm1} and (\ref{eq4}), applied to the convex function $f(x)=x^{\beta}, \
\beta>1$ with $c\le y_i\le d$, we have

$$
0\le \sum_1^n p_i y_i^{\beta}-(\sum_1^n p_i y_i)^{\beta}\le
c(x^{\beta})[c^{\beta}+d^{\beta}-2(\frac{c+d}{2})^{\beta}].
$$

\bigskip

The change of variable $y_i=x_i^{1/\beta}$ gives $a:=c^{\beta}\le
x_i\le d^{\beta}:=b$ and

$$
0\le \sum_1^n p_i x_i-(\sum_1^n p_i x_i^{1/\beta})^{\beta}\le
c(x^{\beta})[ a+b-2(\frac{a^{1/\beta}+b^{1/\beta}}{2})^{\beta}].
$$

\bigskip

Finally, the change of variable $\beta=1/\alpha, \ 0<\alpha<1$,
gives the result.

\bigskip

The second part proof goes analogously, treating the convex
function $f(x)=-x^{\beta}, \ 0<\beta<1$.
\end{proof}

\bigskip

 {\bf A converse of Ky Fan inequality} \ The most celebrated
counterpart of $\mathcal{A}-\mathcal{G}$ inequality is the
inequality of Ky Fan which says that

\begin{equation}
\frac{\sum_1^n w_i x_i}{\sum_1^n w_i (1-x_i)}\ge \frac{\prod_1^n
x_i^{w_i}}{\prod_1^n (1-x_i)^{w_i}}
\end{equation}

whenever $x_i\in (0, 1/2]$.

\bigskip

A converse of Ky Fan inequality is given in \cite{sim}.

\bigskip

\begin{theorem}
If $0<a\le x_i\le b\le 1/2$, then
\begin{equation}
\frac{\sum_1^n w_i x_i}{\sum_1^n w_i (1-x_i)}\le
S(a,b)\frac{\prod_1^n x_i^{w_i}}{\prod_1^n (1-x_i)^{w_i}},
\end{equation}
where

$$S(a, b)=\frac{(1-a)(1-b)(a+b)^2}{ab (2-a-b)^2}$$.
\end{theorem}

\bigskip

A two-sided improvement of this inequality is obtained by an
application of Theorem \ref{thm3}.

\bigskip

\begin{theorem} \ For $0<a\le x_i\le b\le 1/2$, we have

$$
\exp\Bigl(\frac{1/2-b}{(b(1-b))^2}[\sum w_i x_i^2-(\sum w_i
x_i)^2]\Bigr)\frac{\prod_1^n x_i^{w_i}}{\prod_1^n (1-x_i)^{w_i}}
$$
$$
\le \frac{\sum_1^n w_i x_i}{\sum_1^n w_i (1-x_i)}\le
$$
$$
\exp\Bigl(\frac{1/2-a}{(a(1-a))^2}[\sum w_i x_i^2-(\sum w_i
x_i)^2]\Bigr)\frac{\prod_1^n x_i^{w_i}}{\prod_1^n (1-x_i)^{w_i}}.
$$
\end{theorem}

\begin{proof} \ Let $f(x)=\log(\frac{1-x}{x})$. Since
$f''(x)=\frac{1-2x}{(x(1-x))^2}$ and this function is decreasing
on $E=(0,1/2]$, we found that $m_f(E)=\frac{1-2b}{(b(1-b))^2}, \
M_f(E)=\frac{1-2a}{(a(1-a))^2}$.

\bigskip

Therefore, applying Theorem \ref{thm3} we get

$$
\frac{1}{2}m_f(E)J_n(\bold p,\bold x;x^2)\le\sum
w_i\log\Bigl(\frac{1-x_i}{x_i}\Bigr)-\log\Bigl(\frac{1-\sum
w_ix_i}{\sum w_i x_i}\Bigr)
$$
$$
=\log\Bigl(\frac{\sum w_ix_i}{\sum w_i(1-
x_i)}\Bigr)-\log\Bigl(\frac{\prod_1^n x_i^{w_i}}{\prod_1^n
(1-x_i)^{w_i}}\Bigr)\le\frac{1}{2}M_f(E)J_n(\bold p,\bold x;x^2),
$$

and the proof follows.

\end{proof}

\bigskip

It is of interest to find a form of Ky Fan inequality for $\bold
x\in (0,1)$. We shall give now two results of this kind in the
special case $\bold x\in E:=[a, 1-a], 0<a<1/2$.

\bigskip

\begin{theorem}\label{thm2}
If $\bold x\in E:=[a, 1-a], 0<a<1/2$, then

\begin{equation}
\frac{1}{T_n(a; \bold w, \bold x)}\frac{\prod_1^n
x_i^{w_i}}{\prod_1^n (1-x_i)^{w_i}}\le\frac{\sum_1^n w_i
x_i}{\sum_1^n w_i (1-x_i)}\le T_n(a; \ \bold w, \bold
x)\frac{\prod_1^n x_i^{w_i}}{\prod_1^n (1-x_i)^{w_i}},
\end{equation}

where

$$T_n(a; \ \bold w, \bold x)=\exp\Bigl[\frac{1-2a}{2(a(1-a))^2}J_n(\bold w,\bold x;x^2)\Bigr]$$.
\end{theorem}

\bigskip

\begin{proof} \ Analogously to the previous reason, for $f(x)=\log(\frac{1-x}{x})$ we have
$M_f(E)=\frac{1-2a}{(a(1-a))^2}=-m_f(E)$ and the proof is obtained
by Theorem \ref{thm3}. Note that the function $f$ is neither
convex nor concave in this case.
\end{proof}

\bigskip

\begin{corollary} \ A weaker but more explicit variant of the
above assertion is given in the next

\begin{theorem}
If $\bold x\in E:=[a, 1-a], 0<a<1/2$, then

$$
\exp\Bigl[\frac{-(1-2a)^3}{8(a(1-a))^2}\Bigr]\frac{\prod_1^n
x_i^{w_i}}{\prod_1^n (1-x_i)^{w_i}}\le\frac{\sum_1^n w_i
x_i}{\sum_1^n w_i (1-x_i)}\le
\exp\Bigl[\frac{(1-2a)^3}{8(a(1-a))^2}\Bigr]\frac{\prod_1^n
x_i^{w_i}}{\prod_1^n (1-x_i)^{w_i}}.
$$
\end{theorem}

\end{corollary}

\bigskip

\begin{proof} \ Since $c(x^2)=1/2$, we obtain

$$
J_n(\bold w,\bold
x;x^2)\le\frac{1}{2}[a^2+b^2-2(\frac{a+b}{2})^2]=\frac{1}{4}(b-a)^2=\frac{1}{4}(1-2a)^2,
$$

and the result follows from Theorem \ref{thm2}.

\end{proof}

\vskip 1cm

{\bf Applications in Probability Theory} \ The Jensen's inequality
has a great influence in Probability and Statistics. Here are some
basic definitions.

\bigskip

If the generator of random variable $X$ is discrete with
probability mass function $x_1\to p_1, x_2\to p_2,..., x_n\to
p_n$, then the {\it expected value} $EX$ is defined as

$$
EX:=\sum_1^n p_ix_i,
$$

and the {\it variance} $Var(X)$ is

$$
Var(X):=\sum_1^n p_ix_i^2 -(\sum_1^n p_ix_i)^2=E(X^2)-(EX)^2 =
E(X-EX)^2.
$$

\bigskip

Also, the {\it moment of $s$-th order} is defined by

$$
EX^s:=\sum_1^n p_ix_i^s, \ s>0.
$$

\bigskip

{\bf Jensen's moment inequality} says that

$$
EX^s\ge (EX)^s, \ s>1;
$$

and

$$
EX^s\le (EX)^s, \ 0<s<1.
$$

\bigskip

These inequalities follows from the Jensen's inequality applied to
the convex functions $f(x)=-x^s, \ 0<s<1$ and $f(x)=x^s, \ s>1$.
For example $Var(X)\ge 0$.

\bigskip

Our task in the sequel is to improve Jensen's moment inequality by
an application of the results from this paper.

\bigskip

\begin{theorem}\label{thm8} \ For $a\le X\le b$, we have

\begin{equation}\label{eq8}
\frac{1}{2}s(s-1)a^{s-2} Var(X)\le E(X^s)-(EX)^s\le
\frac{1}{2}s(s-1)b^{s-2}Var(X), \ s>2;
\end{equation}

\begin{equation}\label{eq9}
\frac{1}{2}s(s-1)b^{s-2} Var(X)\le E(X^s)-(EX)^s\le
\frac{1}{2}s(s-1)a^{s-2}Var(X), \ 1<s<2;
\end{equation}

\begin{equation}\label{eq10}
\frac{1}{2} s(1-s)b^{s-2} Var(X)\le (EX)^s-E(X^s)\le
\frac{1}{2}s(1-s)a^{s-2}Var(X), \ 0<s<1.
\end{equation}
\end{theorem}

\bigskip

\begin{proof} \ The proof follows by an application of
Theorem \ref{thm3}.

\end{proof}

\bigskip

\begin{theorem}\label{thm9} \ For $a\le X\le b$, we have

\bigskip

\begin{equation}\label{eq11}
0\le (EX)^s-E(X^s)\le
(s-1)s^{s/(1-s)}/(1-2^{1-s})[a^s+b^s-2(\frac{a+b}{2})^s], \ s>1;
\end{equation}

\begin{equation}\label{eq12}
0\le E(X^s)-(EX)^s\le
(1-s)s^{s/(1-s)}/(2^{1-s}-1)[2(\frac{a+b}{2})^s-(a^s+b^s)], \
0<s<1.
\end{equation}
\end{theorem}

\bigskip

\begin{proof} \ Applying (\ref{eq4}) and the result from Theorem
\ref{thm7}, we obtain the proof.

\end{proof}

\bigskip

\begin{remark} \ Comparison of Theorem \ref{thm8} and Theorem
\ref{thm9} is interesting. Although the left-hand side of Theorem
\ref{thm8} is evidently better than the left-hand side of Theorem
\ref{thm9}, what can be said about their right-hand sides?
\end{remark}

\vskip 1cm

 \section{Conclusion}

 The celebrated Jensen's inequality for convex functions is
 applicable in many parts of Analysis, Probability and Statistics,
 Information Theory etc. Some important inequalities such as
 Cauchy's inequality, H$\ddot{o}$lder's inequality, Minkowski's
 inequality, Ky Fan inequality and Jensen-Mercer inequality are just special cases of Jensen;s
 inequality.

 In this article we give several improvements and reverses of
 Jensen's and Jensen-Mercer inequalities. We also consider the
 form of these inequalities for twice differentiable functions
 which are not necessarily convex/concave on a given closed interval.

 Finally, we demonstrate some applications of our results in
 Theory of Means and Probability Theory.

\vskip 1cm

{\bf Acknowledgement} The author is grateful to the referees for
 their valuable comments.

\end{document}